\documentclass[a4paper,10pt]{amsart}
\usepackage[matrix,arrow,tips,curve]{xy}
\usepackage[english]{babel} 
\usepackage[leqno]{amsmath}
\usepackage{amssymb,amsthm}
\usepackage{amscd}
\usepackage{enumerate}

\newtheorem{thm}{Theorem}[]
\newtheorem{lemma}[thm]{Lemma}

\newtheorem{prop}[thm]{Proposition}

\newtheorem{claim}[thm]{Claim}
\newtheorem{cor}[thm]{Corollary}

\theoremstyle{remark}
\newtheorem{remark}[thm]{Remark}
\theoremstyle{definition}
\newtheorem{defi}[thm]{Definition}

\newcommand{\la}{\longrightarrow}

\newcommand{\da}{\dashrightarrow}

\newcommand{\ov}{\overline}

\newcommand{\Div}{\operatorname{Div}}

\newcommand{\im}{\operatorname{Im}}

\newcommand{\Pic}{\operatorname{Pic}}
\newcommand{\supp}{\operatorname{Supp}}
\newcommand{\mdeg}{\operatorname{{\underline{de}g}}}

\def\O{\mathcal O}

\newcommand{\Z}{\mathbb{Z}}

\def\me{\underline{e}}
\def\md{\underline{d}}

\newcommand{\pr}[1]{\mathbb{P}^{#1}}

\newcommand{\Mgb}{\ov{M}_g}
\newcommand{\Mgrd}{M^r_{g,d}}
\newcommand{\Mgrdb}{\ov{\Mgrd}}
\newcommand{\Hgb}{\ov{H}_g}

\def\Pdgb{\overline{P}_{d,\,g}}

\newcommand{\dfib}{\dot{X}}

\def\PXb{\overline{P^d_X}}

\def\PXdb{\overline{P^2_X}}

\def\Wmd{W_{\md}(X)}
\def\Wmdr{W_{\md}^r(X)}

\def\BN{\rho^r_d(g)}

\def\Amd{A_{\md}(X)}
\def\Xmd{X^{\md}}

\def\amd{\alpha^{\md}_X}

\def\FM{F_{M}(X)}

\newcommand{\picX}[1]{\Pic^{#1}X}

\def\mgbst{\overline{\mathcal{M}}_g}

\newcommand{\YS}{Y_S}

\newcommand{\hXS}{{\widehat{X}}_S}
\newcommand{\hLS}{{\widehat{L}}_S}
\newcommand{\hL}{{\widehat{L}}}

\newcommand{\hmd}{{\widehat{\md}}}

\newcommand{\sing}{X_{\text{sing}}}

\def\BXd{B_d  (X)}
\def\Bgd{B_d  (g)}
\def\BXds{B_d^*  (X)}
\def\Bgds{B_d^*  (g)}
\def\BYs{B_{d-e}^*  (Y_{S^e})}
\def\BX1{B_1  (X)}
\def\BXt{B_3  (X)}

\newcommand{\Cl}{\operatorname{Cliff}}
\newcommand{\Bg}{B_g}
\newcommand{\Bgrd}{B_{g,d}^r}
\newcommand{\Bgrmd}{B_{g,{\underline d}}^r}
\newcommand{\Bgtmd}{B_{g,{\underline d}}^2}

\newcommand{\Bgtd}{B_{g,d}^2}

\newcommand{\Mn}{M_{0,n}}

\newcommand{\Mapr}{M_{0}(\pr{r},d)}

\newcommand{\Wb}{\ov{W^r_{d,X}}}

\newcommand{\Wudb}{\ov{W^1_{2,X}}}

 \address{Dipartimento di Matematica \\ Universit\`a Roma 3\\ Largo S.L. Murialdo 1, 00146 Roma  Italy}
 \email{caporaso@mat.uniroma3.it}

\begin{document}
\begin{center}
{\bf\large Brill-Noether theory of binary curves }

{Lucia Caporaso}
\end{center}

\noindent
{\bf  Abstract.} 
The theorems of Riemann, Clifford  and Martens are proved for every line bundle
parametrized by the compactified Jacobian of every binary curve.
The Clifford index is used to characterize hyperelliptic and trigonal binary curves.
The  
Brill-Noether theorem  for $r\leq 2$  is proved for a general binary curve.

\section{Introduction}

The purpose of this paper  is to
contribute to the
Brill-Noether theory of stable curves, about which very little is known. 
We work over an algebraically closed field, and consider the
 compactified  universal Picard variety,
$\Pdgb \to \Mgb$, parametrizing degree-$d$ balanced line bundles on semistable curves of genus $g$
(or, which is equivalent, semistable torsion-free sheaves of rank one on stable curves).
The moduli properties of $\Pdgb$ are nowadays  quite well understood, 
both from the scheme theoretic point of view and the stack theoretic one; moreover it has several equivalent geometric descriptions
\cite{alex},   \cite{melo}, \cite{cner}.
In this paper, the Brill-Noether varieties of stable curves are     defined  inside $\Pdgb$.

In older times, lacking a thorough
understanding of how to compactify the Picard functor,
or, later on, in the presence of   different, seemingly unrelated, solutions of this problem,
research    about such topics
 followed  different   approaches. As examples,  
let us recall two  famous constructions,
which have   had   several important applications. 
The first is the theory of admissible covers, due to J. Harris and D. Mumford
\cite{HM},
studying degenerations of   linear series of dimension one.
The second is  the theory of limit linear series, created  by D.  Eisenbud
and J. Harris \cite{EH};
this theory, valid for linear series of any dimension,
 makes no use of compactified Jacobians, and works  
 best for curves of compact type,  whose Jacobian is projective; see also \cite{bruno}, \cite{EM}
and  \cite{oss} for more recent developements. 

The subsequent progress on compactified moduli spaces of line bundles 
 followed different directions.
This   led to the
construction of moduli spaces (the compactified Jacobians, or Picard schemes, mentioned at the beginning) which are
  natural ambient spaces where  studying Brill-Noether type questions.
  
In this field  there are many open problems, some of which
 appear almost intractable, 
owing to   the  combinatorial complexity of  stable curves.
As a consequence, much of the previous work on the subject deals only with certain types of stable curves:
(\cite{EH}, \cite{oss}
dealing with curves of compact type, or \cite{bruno}, 
\cite{EM} dealing with curves with two components).
In the present paper also, only a  certain type of curve is studied: the  so-called ``binary
curves", namely, nodal curves made of two smooth rational components, intersecting at $g+1$ points.
Their moduli scheme   is  irreducible of dimension $2g-4$.

Binary curves   arise naturally in a variety of situations, 
sometimes with a different name, such as ``split curves". Their canonical model
(for non-hyperelliptic ones) is the union of two rational normal curves meeting transversally at $g+1$
points, a remarkable   curve,  useful   as a test case and as a limit case.
Also, canonical  binary curves specialize to rational ribbons, another particularly interesting type of curve.

Although binary curves are reducible,   many   numerical and combinatorial difficulties tremendously
simplify for them. 
Moreover, as they are made of  rational components, moduli spaces of marked rational curves,
and of their maps   to projective spaces,
provide a powerful tool.

We begin the paper with  some
preliminary results about compactified Jacobians and
  Brill-Noether varieties.
Then we proceed to extend some among the fundamental theorems
on which the classical Brill-Noether theory of Riemann surfaces is based:
the theorems of Riemann, of Clifford, of Martens, and of Brill-Noether.
Notice that none of them is known for all stable curves.

The first three of them
are here proved to hold for the line bundles parametrized by the compactified
Picard scheme.
The  analogue of   Riemann's theorem  
 is not difficult; 
see Proposition~\ref{br}.
In Section~\ref{clifsec} we establish   Clifford's   theorem  and 
study the Clifford index
(Theorem~\ref{bcl}), characterizing binary curves
having Clifford index $0$ or $1$
in terms of their gonality. 
We extend    Martens theorem in
Proposition~\ref{martens}.
We never use   that such theorems hold for smooth curves.

While the rest of the paper deals with every binary curve,
Section~\ref{BNsec} focuses on the general one  and   is devoted to the
Brill-Noether theorem (on the dimension of Brill-Noether varieties) for
$r\leq 2$; see Theorem~\ref{BNt}.  The proof is  independent from the    theorem for smooth curves,
which can hence be re-obtained as a consequence.

Finally, a few words about  further developments.
For binary curves, there are   several  appealing  questions remaining, such as a Brill-Noether
theorem for higher $r$.
Another  direction is to consider all stable curves: how do our results generalize?
In both cases the
situation is considerably more complex; in fact, our preliminary investigation
(to appear in  a forthcoming paper)  has shown
that the Clifford theorem 
and the Brill-Noether theorem do fail in some cases.

\section{Set-up}
\label{setup}
\subsection{Binary curves and balanced line bundles}
A  reduced nodal
curve $X$ is called a {\it binary curve} if $X=C_1\cup C_2$ with $C_i\cong \pr{1}$; let $g$ be 
the arithmetic genus of
$X$, then  $g\geq -1$ and
$\#C_1\cap  C_2 =g+1$.

As binary curves are union of smooth rational components, 
certain  moduli spaces 
come naturally into the
picture. 
For any $n\geq 4$ consider $\Mn$, the moduli space of $n$-marked  smooth rational curves.
$\Mn$ is irreducible of dimension $n-3$.

 We denote by $\Mapr$  the moduli space of maps of degree $d\geq 1$ from $\pr{1}$ to
$\pr{r}$. 
More generally, for any $n\geq 0$
consider 
the moduli space $M_{0,n}(\pr{r},d)$ of degree $d$ maps    from   $n$-marked, smooth, rational curves to $\pr{r}$.
It is irreducible of dimension 
\begin{equation}
\label{maps}
\dim M_{0,n}(\pr{r},d)=\dim \Mapr+n=(r+1)d+r-3+n.
\end{equation}

\begin{lemma}
\label{Bg}
Let $\Bg\subset \Mgb$ be the locus of binary curves of genus $g\geq 2$.
 Then
$\Bg$ is irreducible of dimension $2g-4$. 
\end{lemma}
\begin{proof} 
There is a surjective morphism, having finite fibers,
\begin{equation}
\label{proj}
\gamma_g:M_{0, g+1}\times M_{0, g+1} \la \Bg    
\end{equation}
mapping  
$\Bigr( (C_1;p_1,\ldots ,p_{g-2},0,1,\infty),(C_2;q_1,\ldots ,q_{g-2},0,1,\infty)\Bigl)$
to the binary curve obtained by gluing $p_i$ with $q_i$ and $0,1,\infty \in C_1$
with $0,1,\infty \in C_2$. As $M_{0, g+1}$ is irreducible of dimension $g-2$, the Lemma follows.
\end{proof}

The description of the compactified Picard scheme of a binary curve
(see Section~\ref{compP} below) is based on  
Definition~\ref{bal},
 a special case of (for example)  4.6 in \cite{cner}.

\begin{defi}\label{bal}
Let $X$ be a binary curve of genus $g\geq -1$.
A multidegree $\md=(d_1,d_2)$ with  $d=|\md|=d_1+d_2$  is {\it balanced} on $X$ if,  for either $ i\in \{1,2\}$,
\begin{equation}
\label{BIb}
m(d,g):=\frac{d-g-1}{2}\leq d_i \leq \frac{d+g+1}{2}=:M(d,g).
\end{equation}
We say that   $L\in \picX{d}$
is  balanced if  $\mdeg L $ is balanced on $X$.
We say that $\md$, or $L$, is  {\it strictly balanced}
if (\ref{BIb}) holds with strict inequalities.
We denote
\begin{equation}
\label{bals}
\BXd=\{\md: |\md|=d, \  \md \text{ balanced }\}\supset \BXds=\{\md \text{ strictly balanced}\}.
\end{equation}
\end{defi}
Clearly $\BXd$ and $\BXds$ depend only on $g$, so we shall sometimes write
\begin{equation}
\label{balg}
\Bgd:=\BXd,\  \  \   \Bgds:=\BXds.
\end{equation}
\begin{remark}\label{balr}
The following  facts will be used several times.
\begin{enumerate}[(a)] 
\item\label{balra}
For every $d$:  $\BXds \neq \emptyset$ if $g\geq 1$, and $\BXd \neq \emptyset$
if $g\geq 0$.
\item\label{balrb}
If $g=-1$, then $\BXd \neq \emptyset \iff m(d,g)\in \Z$.
\item\label{balrc}
$\md\in \BXd \iff
 d_i\geq m(d,g)$,\   $\forall i=1,2\iff
 d_i\leq M(d,g)$,\   $\forall i=1,2.$
\item
\label{balrd}
$\md$ is balanced $\iff$ $\md+n\mdeg\omega_X$ is balanced.
\end{enumerate}
\end{remark}
\begin{remark}
\label{mM}
Let $d$ and $g\geq -1$ be  integers.
Then one easily checks the following.
\begin{enumerate}[(A)]
\item
\label{mMa} $m(d,g)= m(d-1,g-1)$ and $M(d,g)= M(d-1,g-1)+1$.
\item
\label{mMc}
$m(d,g)> m(d,g+n)$ for every $n\geq 1$.
\item
\label{mMd}
$M(d,g)< M(d,g+n)$ for every $n\geq 1$.
\item
\label{ugb}
$B_d(g)\subset B_d(g+n)$ for any $n\geq 0$.
\end{enumerate}
\end{remark}

As it is well known, there are two common (equivalent) ways of describing the geometric objects 
parametrized by the compactified Jacobian:
via torsion-free sheaves or via line bundles; we choose the second one, introduced in \cite{caporaso}.
In order to describe it, we introduce some terminology.
Let   $X$ be a  nodal  curve  and   $S$   a set of nodes of $X$.
By ``the normalization of $X$ at $S$" we mean the local desingularization (or normalization)
of $X$ at every node in $S$.  We say that a nodal curve $\hXS$ is the ``blow -up"  of $X$ at
$S$ if there exists $\pi: \hXS\to X$ such that
$\pi^{-1}(n_i)=E_i\simeq\pr{1}$ for
any $n_i\in S$, and $\pi: \hXS\smallsetminus\cup_iE_i\to 
X\smallsetminus S$
is an isomorphism. Thus $\ov{\hXS\smallsetminus\cup_iE_i}$ is the normalization of $X$ at $S$. 

The boundary points of the compactified Jacobian  parametrize   balanced line bundles on (strictly) semistable curves;
a balanced line bundle is defined to be one whose multidegree is balanced.
To define this for strictly semistable curves we introduce some notation that will be used throughout the paper.
Let $X$ be a binary curve and
  $S\subset \sing$  be a set of  nodes of $X$, set $e=\#S$; we shall sometimes  write $S=S^e$.
We denote $\hXS$ the blow-up of   $X$ at $S$. We call $E_1,\ldots, E_e$ the exceptional
components of
$\hXS$, and
$Y_S$ their complementary curve (the normalization of $X$ at $S$).
 $Y_S$ is  a  binary curve of genus $g-e$, and  
$$ \hXS=Y_S \cup\cup_{i=1}^eE_i=C_1\cup C_2\cup E_1\cup\ldots \cup E_e.$$ 
We will write a multidegree $\hmd=(d_1,d_2,d_3,\ldots, d_{2+e})$  on $\hXS$
using the convention that for $i=1,2$ we have $d_i=\hmd_{C_i}$, and for $i\geq 3$ we have $d_i=\hmd_{E_i}$.
We also write $\hmd_{Y_S}=(d_1,d_2)$, so that $|\hmd_{Y_S}|=d-\sum_{i=3}^ed_i$.
\begin{defi}
\label{urb}
A multidegree  $\hmd$  on $\hXS$, with $|\hmd|=d$, is   balanced   if (1) and (2) hold:

\noindent
(1) $d_i=1,$ $\forall i=3,\ldots,e$  
(i.e. if $\hmd_{E_i}=1$ for every $E_i$);
 
\noindent
(2) $\hmd_{Y_S}$ is  balanced on $Y_S$  (i.e. if
$\hmd_{Y_S}\in B_{d-e}(Y_S)$).

\noindent
$\hmd$   is  called strictly  balanced if $\hmd_{Y_S}$ is strictly  balanced on $Y_S$.
\end{defi}
We denote  $B_d(\hXS)$ and $B_d^*(\hXS)$ the set of balanced and strictly balanced multidegrees on $\hXS$.
As we said, $\hL\in \Pic^{\hmd}\hXS$ is called balanced if $\hmd$ is balanced.
Two balanced line bundles $\hL',\hL\in \Pic^{\hmd}\hXS$  are defined to be {\it equivalent}
if their restrictions to $\YS$ are isomorphic.

\subsection{The compactified Picard scheme of binary curves.}
\label{compP}
Let $X$ be a stable binary curve of genus $g\geq 2$, and   $d$  a fixed integer.
We shall now describe its compactified degree-$d$ Picard variety
$\PXb$. 
As $d$ varies, the structure of $\PXb$
varies between two different types, according to whether or not $m(d,g)$
is an integer. The terminology we will use reflects the relation with N\'eron models; see \cite{cner} and \cite{melo}.

\noindent{\bf  N-type: } $m(d,g)\not\in \Z$. $X$ is said to be $d$-{\it general}, and $\PXb$ of {\it
N\'eron type}.

In this case      every point of  $\PXb$ corresponds to an equivalence
class of balanced 
line bundles.
We have a natural isomorphism
\begin{equation}
\label{Pstr1}
\PXb\cong \coprod _{\md \in \BXd}\picX{\md}\coprod_{e=1}^g\Bigr(
\coprod_{\stackrel{S^e\subset \sing}{\#S^e=e}}\coprod_{\  \  \md^e\in B_{d-e}(Y_{S^e})}
\Pic^{\md^e}Y_{S^e}\Bigl).
\end{equation}
Note that $\BYs=B_{d-e}(Y_{S^e})$ for every $\emptyset\subseteq S^e\subset \sing$.

\noindent
{\bf  D-type: } $m(d,g)\in \Z$. Now $\PXb$ is called of  {\it Degeneration
type}.

In this case there exist balanced multidegrees that are not strictly balanced. 
More precisely, for every partial normalization $Y_{S^e}$ of $X$, $e\geq 0$,
there exists a unique such multidegree,
namely   
$(m(d,g), M(d,g)-e)\in B_{d-e}(Y_{S^e})$ (cf. Lemma~\ref{mM}).
All   line bundles having these multidegrees are identified to a unique point $\ell_0\in \PXb$.
Of course, to $\ell_0$ there corresponds    a unique closed orbit; indeed 
there exists a unique balanced line bundle on a unique curve parametrized by $\ell_0$,
namely the line bundle $(\O_{C_1}(m(d,g)),\O_{C_2}(m(d,g)))$ on the normalization of $X$
(the disjoint union of two copies of $\pr{1}$).
We have a description analogous to (\ref{Pstr1})
\begin{equation}
\label{Pstr2}
\PXb\smallsetminus\{\ell_0\}\cong \coprod _{\md \in \BXds}\picX{\md}\coprod_{e=1}^{g-1}\Bigr(
\coprod_{\stackrel{S^e\subset \sing}{\#S^e=e}}\coprod_{\  \  \md^e\in \BYs}
\Pic^{\md^e}Y_{S^e}\Bigl).
\end{equation}
Note that if $e=g$ then $\BYs$ is empty.

\

For any $S^e\subset \sing$ and any $\md^e\in B_{d-e}(Y_S)$ 
we shall denote $P_{S^e}^{\md^e}\subset \PXb$ the stratum isomorphic to $\Pic^{\md^e}Y_{S^e}$.
Also, for a fixed $S\subset \sing$  
we denote $P_S$ the union of all strata $P_{S}^{\md}$ ad $\md$ varies,
 omitting ``$e$" from the notation, for simplicity.
Note that all the strata above are tori:
$
P_{S^e}^{\md^e}\cong (k^*)^{g-e}.
$
Moreover,  
\begin{equation}
\label{filt}
\overline{P_{S}^{\md}}\supset P_{S'}^{\md'} \iff S\subset S' \text{ and } \md\geq \md'
\end{equation}
where $\md'\in B_{d-e'}(Y_{S'})$, and $\md\geq \md'$ means $d_i\geq d'_i$, $i=1,2$.

\subsection{Brill-Noether varieties.}
\label{notBN}
Given $\md$ and $r$ we denote
\begin{equation}
\label{Wmd}
 W_{\md}^r(X):=\{L\in \picX{\md}: h^0(L)\geq  r+1\}
\end{equation}
if $r=0$ we usually omit $r$: $ W_{\md}^0(X)=\Wmd=\{L\in \picX{\md}: h^0(L)\neq 0\}$.
$W_{\md}^r(X)$ is endowed with a natural scheme structure,  obtained either as for smooth curves
(\cite{ACGH}), or using the GIT construction of $\PXb$. We omit the details as this is irrelevant for our purposes.
For any $r$ and $d$, we denote
\begin{equation}
\label{Bgrd}
\Bgrd=\{X\in \Bg : \exists\md \in \BXd: \  W^r_{\md}(X)\neq \emptyset\},
\end{equation}
and for any $\md\in \Bgd$
\begin{equation}
\label{Bgrmd}
\Bgrmd=\{X\in \Bg :   \  W^r_{\md}(X)\neq \emptyset\}.
\end{equation}
By the   above description,  
every point $\lambda$ of $ \PXb$,
 belongs to a stratum $P_S$, for some $S\subset \sing$. So 
$\lambda$ determines   a unique strictly balanced line bundle $M_S$,
of degree $d-\#S$,  on a unique  curve
$Y_S$,   the  normalization of $X$ at $S$.
Viewing the isomorhisms of (\ref{Pstr1}) and (\ref{Pstr2}) as identifications, we shall often denote
the points of $\PXb$ as follows
\begin{equation}
\label{notP}
 [M,S]\in \PXb,\  \   S\subset \sing, \  \  M\in \Pic^{d-\#S}Y_S,\  
\end{equation}
where $M$ is strictly balanced.
On the other hand, a point
   of $\PXb$    parametrizes a pair
$(\hXS, [\hL] )$,
where  
  $\hXS=Y_S\cup_{i=1}^{\#S}E_i$ is the  blow-up of $X$  at $S$, and 
$[\hL]$ is an equivalence class of strictly balanced line bundles on $\Pic^d\hXS$,
all having restriction $M$ on $Y_S$.
So, we will also  denote simply by $[\hL]$ a point of  $\PXb$.

With the above notations, one easily sees (cf. Lemma 4.2.5 in \cite{Ctheta})
\begin{equation}
\label{key}
h^0(Y_S, M)=h^0(\hXS, \hL '),\  \   \forall \hL '\in [\hL_S].
\end{equation}
Now, we define
\begin{equation}
\label{Wb}
\Wb=\{[M,S]\in \PXb: h^0(Y_S, M)> r\}=\{[\hL]\in \PXb: h^0(\hXS, \hL)> r\}.
\end{equation}

We denote by $\Mgrd\subset M_g$ the locus of smooth curves $C$ such that $W^r_d(C)\neq \emptyset$,
and by $\Mgrdb\subset \Mgb$ its closure in $\Mgb$.
\begin{prop}
\label{Wp} Let $r\geq 0$, $g\geq 2$ and $d\leq r+g-1$.
\begin{enumerate}[(i)]
\item
\label{}
There is a natural isomorphism
\begin{equation}
\label{Wstr}
\Wb\cong   \coprod _{\md \in \BXds}\Wmdr\coprod_{e=1}^g\Bigr(
\coprod_{\stackrel{S^e\subset \sing}{\#S^e=e}}\coprod_{\  \  \md^e\in B_{d-e}^*(Y_{S^e})}
W^r_{\md^e}(Y_{S^e})\Bigl).
\end{equation}
\item
\label{Wp2}
Denote by $W_{\md^e,S^e}\subset \Wb$ the stratum isomorphic to $W^r_{\md^e}(Y_{S^e})$ under (\ref{Wstr}).
If  $\overline{W_{\md^e,S^e}}\supset W_{\md^{e'},S^{e'}}$ then $S^e\subset S^{e'}$ and $\md^e \geq \md^{e'}$.
\item
\label{Wp3} If $X\in \Mgrdb\  $  then $\Wb\neq \emptyset$.
\end{enumerate}
\end{prop}
\begin{proof}
We earlier gave  a description of $\PXb$ by a natural isomorphism analogous to (\ref{Wstr}).
We explained that there are two possibilities, (\ref{Pstr1}) and (\ref{Pstr2}),
according to whether $m(d,g)$ is an integer or not. If $m(d,g)\not\in \Z$, then 
(\ref{Pstr1}) holds and  $\BYs=B_{d-e}(Y_{S^e})$ for every $e$ and $S$.
Therefore (\ref{Wstr}) follows immediately from (\ref{Pstr1}).

Suppose $m(d,g)\in \Z$, and consider the line bundle
$$
M:=\Bigl(\O_{C_1}(m(d,g)) ,\O_{C_2}(m(d,g))\Bigr)\in \Pic (C_1\coprod C_2)
$$
corresponding to the point $\ell_0\in \PXb$.
Now, as $d\leq r+g-1$, we have
$$
m(d,g)=\frac{d-g-1}{2}\leq \frac{r+g-1-g-1}{2}=\frac{r}{2}-1.
$$
Therefore 
$ h^0(M)=2h^0(\pr{1}, \O_{\pr{1}}(m(d,g)))  =r
$ 
hence $\ell_0\not\in\Wb$.
This implies that (\ref{Wstr}) follows   from (\ref{Pstr2}), as in the previous case.

Part (\ref{Wp2}) follows   from the previous one and from (\ref{filt}).

Now part (\ref{Wp3}).
Let $X\in \Mgrdb$.
Then there exists a family of smooth curves specializing to $X$ such that the general fiber, $C$, of the family has
a non empty $W^r_d(C)$. Up to replacing the family by some base change, we may assume that the family has a
section. This enables us to apply a construction of E. Arbarello and M. Cornalba (see Section 2 of  \cite{AC1})
yielding that the $W^r_d(C)$ form a family contained in the relative Picard scheme.
Therefore, as $C$ specializes to $X$, $W^r_d(C)$ specializes to some non-empty subset $W_0$ of $\PXb$.
By uppersemicontinuity of $h^0$, $W_0$ lies in $\Wb$, which is thus non empty.
\end{proof}

For every $d\geq 1$ and $\md=(d_1,d_2)$,
denote $\Xmd:=C_1^{d_1}\times C_2^{d_2}$. Consider the Abel map of multidegree $\md$
\begin{equation}\label{}
\amd:C_1^{d_1}\times C_2^{d_2}\da \Wmd; \  \  \  (p_1,\ldots,p_d)\mapsto \O_X(\sum_{i=1}^d p_i). 
\end{equation}
$\amd$ is regular away from the points lying over $C_1\cap C_2$.
We denote $\Amd\subset \Wmd$ the closure of the image of $\amd$. It is clear that $\Amd$ 
is irreducible.
\begin{lemma}
\label{abel}
Let $1\leq d\leq g$ and $\md \in \BXd$. Then
 $h^0(X,L)=1$  for the general $L\in \Amd$, and $\dim \Amd =d$.
\end{lemma}
\begin{proof} We have $\dim \Amd \leq d$, of course.
The fiber of $\amd$ over a general $L\in \Amd$ has dimension $h^0(L)-1$, hence
it suffices to prove that
$h^0(L)\leq 1$ for some $L\in \Amd$. 

Pick $S\subset \sing$ such that $\#S=d$. As $d<g+1=\#\sing$ we can consider
the normalization of $X$  at $S$, $Y_S\to X$, and the curve $\hXS$, the blow-up of $X$ at $S$.

Consider $M_S\in \Pic^{(0,0)}Y_S$, note that, since $Y_S$ is connected, $h^0(Y_S,M_S)\leq 1$,
and equality holds if and only if $M_S=\O_{Y_S}$.
Therefore, as   already observed in (\ref{key}), for every balanced line bundle $\hL$ on $\hXS$ restricting to $M_S$
on $Y_S$, we have
\begin{equation}
\label{leq1}
h^0(\hXS, \hLS)=h^0(Y_S,M_S)\leq 1 
\end{equation}
with equality if and only if $M_S=\O_{Y_S}$
(by Corollary 2.2.5 of \cite{Ctheta}).
Fix $M_S=\O_{Y_S}$ and $\hLS$ as above. 
The point of $\PXb$ parametrizing $\hLS$ is in the closure of
$\Amd$. Indeed, we can simultaneously specialize $d$ distinct nonsingular points 
of $X$ to the $d$ nodes of $S$.
By (\ref{leq1})  we get $h^0(X,L)\leq h^0(\hXS, \hLS)\leq 1 $ for $L$ general in $\Amd$, as wanted.
\end{proof}

Let $\nu:Y\to X$ be  the normalization of $X$ at $s$  nodes, $n_1,\ldots,n_s$,
and $\nu^{-1}(n_s)=\{p_s,q_s\}$. In symbols:
 \begin{equation}\label{pm}
\nu:Y\la
X=Y/_{\{p_i=q_i, \  i=1\ldots s\}}.
\end{equation}
Let $M$ be a line bundle on $Y$ such that $h^0(Y,M)\neq \emptyset$.
Denote by $\FM$ the fiber over $M$ of  
$\nu^*:\Pic X\to \Pic Y$, i.e.
$
\FM :=\{L\in \Pic X:\nu^*L=M\}.
$
We ask under what conditions there exists $L\in \FM$ such that  
$h^0(X,L)=h^0(Y,M)$.
We introduce the following terminology.
\begin{defi}
\label{np}
Let $p,q$ be nonsingular points of a curve $Y$; pick $M\in \Pic Y$.
We say that $p$ and $q$ are {\it equivalent, or   neutral, with respect to $M$},
and write
$ 
p\sim_M q, 
$ 
if $
h^0(Y,M-p)=h^0(Y,M-q)=h^0(Y,M-p-q).
$ 
\end{defi}

\noindent
The following  is a straightforward consequence  of Lemmas 2.2.3 and  2.2.4 in \cite{Ctheta}.

\begin{lemma}
\label{dsl} Let $Y\to X=Y/_{\{p_1=q_1,\ldots ,p_s=q_s\}}$;  
pick $M\in \Pic Y$  with  $h^0(Y,M)\neq 0$.  
There exists $L\in \FM$ such that
$h^0(X,L)=h^0(Y,M)$ if and only if 
$ 
p_i\sim_M q_i 
$ 
for every $i=1,\ldots, s$.

Such an $L$   is unique (if it exists) if 
   $p_i$ and $q_i$ are not base points for $M$ for all $i$.
\end{lemma}

This implies the following  useful result.
\begin{lemma}
\label{e} 
Let $\md=(d_1,d_2)$ be a  multidegree on a binary curve $X$ of genus $g$; assume $d_2\geq d_1\geq -1$.
Then for every $L\in \picX{\md}$
\begin{equation}
\label{eeq}
h^0(X,L)\leq d_1+d_2+1-\min \{d_2,g\}.
\end{equation}
\begin{enumerate}[(i)]
\item
\label{geq}
If $d_2\geq g$, equality holds 
for every $L\in \picX{\md}$.
\item
\label{!}
If $d_2< g$, equality holds 
for at most one $L\in \picX{\md}$.
\end{enumerate}
\end{lemma}
\begin{proof}
Set $d=d_1+d_2$. 
If $g=-1$ then $\min\{d_2,g\}=-1$,  hence (as $d_i\geq -1$)
$$
h^0(X,L)=h^0(C_1,L_1)+h^0(C_2,L_2)=d_1+1+d_2+1=d+1-\min\{d_2,g\}.
$$
We can assume    $g\geq 0$, i.e.  $X$ is connected. For every $0\leq e\leq \min \{d_2,g\}$,
denote 
$$
X_e:=\frac{C_1\coprod C_2}{(p_i=q_i,\   i=1,\ldots, e)} \  \stackrel{\nu_e}{\la}X
$$
so that $\nu_e$ is a normalization at $g+1-e$ nodes.
Set $M_e=\nu_e^*L$; we have, of course, $h^0(X_e,M_e)\geq h^0(X,L)$.

If $e=0$ then $h^0(X_0,M_0)=d_1+d_2+2=d+2$. More generally, we claim that 
$$
h^0(X_e,M_e)=d+2-e
$$ 
 for every $e$. By induction on $e$.
Notice that  $\deg L_2(-\sum_{i=1}^eq_i)\geq 0$,
therefore, as $C_2\cong \pr{1}$,
 there exists a section 
$s_2\in H^0(C_2, L_2(-\sum_{i=1}^eq_i))$
not vanishing at $q_{e+1}$. This implies that $M_e$ has a section vanishing
at $p_{e+1}$ but not at $q_{e+1}$; indeed, just glue $s_2$ to the zero section on $C_1$, which we can do as $s_2$
vanishes at every $p_i$ with $i\leq e$. Therefore
 $p_{e+1}\not\sim_{M_e} q_{e+1}$.
 Lemma~\ref{dsl} now yields
$$
h^0(X_{e+1},M_{e+1})= h^0(X_{e},M_{e})-1=d+2-e-1=d+1-e.
$$
Applying this to $e=\min \{d_2,g\}$
we obtain
$$
h^0(X,L)\leq h^0(X_{e+1},M_{e+1})=d+1-\min \{d_2,g\}.
$$
We have thus shown that (\ref{eeq}) holds, with equality if $d_2\geq g$.

Part (\ref{!}) follows from the uniqueness part in Lemma~\ref{dsl}.
\end{proof}

Using Lemma~\ref{e} we can now extend
Riemann's theorem:
\begin{prop}
\label{br} 
Let $X$ be a  binary curve   of genus $g$, and let $d\geq 2g-1$.

\begin{enumerate}[(i)]
\item
\label{brs}
For every balanced 
$L\in \picX{d}$ we have
$
h^0(L)=d-g+1.
$
\item
\label{brss}
For every $[\hL]\in \PXb$ we have $h^0(\hL)=d-g+1$.
\end{enumerate}

\end{prop}
\begin{proof}
Let $\mdeg L =(d_1,d_2)$ and assume $d_1\leq d_2$. Then 
$
d_2\geq g 
$,
for otherwise
$d_1+d_2\leq 2(g-1)$ which is ruled out, by hypothesis.
As $\mdeg L$ balanced, we  have 
$$
d_i\geq m(d,g)=\frac{d-g-1}{2}\geq \frac{2g-1-g-1}{2}=\frac{g-2}{2}\geq -\frac{3}{2}.
$$
Therefore $d_i\geq -1$ for $i=1,2$, so that
 Lemma~\ref{e} applies.  We obtain that (\ref{eeq}) holds, with equality, as $d_2\geq g$. Hence 
$$
h^0(X,L)=d_1+d_2+1-\min\{d_2,g\}=d+1-g,
$$
as stated in part (\ref{brs}).
Now, to prove  part (\ref{brss}) it suffices to consider $\hL\in \Pic \hXS$ with $\#S=e\geq 1$
(notation as in Subsection~\ref{notBN}).
By (\ref{key}) we have
$$
h^0(\hXS,\hL)=h^0(Y_S,M)=(d-\#S)-(g-\#S)+1=d-g+1
$$
where the second equality follows from   part (\ref{brs})  applied to the binary curve $Y_S$ (of course
$Y_S$ has  genus 
$g-e$, so that $d-e\geq 2g-1-e\geq 2g-1-2e=2g_{Y_S}-1$).
\end{proof}

\begin{prop}
\label{empty} 
Let $\md=(d_1,d_2)$ be a  balanced multidegree on a binary curve $X$.
Assume $d_1\leq d_2$ and set $d=|\md|$. Then
$\Wmdr =\emptyset$ in the following cases.
\begin{enumerate}[(i)]
\item
\label{e<}
$d_1<0$ and $d\leq g+r$.
\item
\label{e>}
$0\leq d_1\leq r-1$ and $d\leq g+r-1$.
\end{enumerate}
\end{prop}
\begin{proof}
We must prove that  $h^0(X,L)\leq r$ for every  $L\in \picX{\md}$.
In case (\ref{e<})
$$
h^0(X,L)=h^0(C_2,L_2(-C_1\cap C_2))=\max\{0, d_2-g\}.
$$
$\md$ is balanced, hence  
$$
d_2-g\leq (d+g+1)/2-g=(d-g+1)/2\leq (r+1)/2.
$$
We obtain $h^0(X,L)\leq  (r+1)/2 <r+1$ and we are done.

In case (\ref{e>}), as $d_1\leq d_2$, we
 have, by Lemma~\ref{e},
$ 
h^0(X,L)\leq d+1-\min\{d_2,g\}.
$ 

If $d_2\leq g$ we obtain
$ 
h^0(X,L)\leq d+1-d_2=d_1+1\leq r-1+1=r
$ 
and we are done.
If $d_2> g$
we have
$ 
h^0(X,L)\leq d+1-g\leq r.
$ 
The proof is complete. 
\end{proof}

\section{Clifford theory}
\label{clifsec}
 \subsection{Clifford's inequality and hyperelliptic binary curves.}
The main result of this Section is Theorem~\ref{bcl},   extending   Clifford's theorem. Its first part, the
Clifford inequality, is  the subsequent Proposition~\ref{cliff}.

\begin{prop}[Clifford's inequality]
\label{cliff}
Let $X$ be a binary curve of genus $g\geq 1$, and let
  $d$ be such that $0\leq d\leq 2g$.

\begin{enumerate}[1.]
\item
\label{cls}
For every $\md \in \BXd$, and every $L\in \picX{\md}$,
we have $h^0(L)\leq d/2+1$.

\noindent
If $d=0$ and $h^0(L)=1$ then $L=\O_X$;
if $d=2g-2$ and $h^0(L)=g$ then $L=\omega_X$.
\item
\label{clss}
For every $[\hL]\in \PXb$ we have $h^0(\hL)\leq d/2+1$.
\end{enumerate}
\end{prop}
\begin{proof}
We may assume $d_1\leq d_2$.
If $d_1<0$ then
$$
h^0(X,L)=h^0(C_2,L_2(-C_1\cap C_2))=d_2-g\leq M(d,g)-g=\frac{d-g+1}{2}
$$
($\md$ is balanced). Now, as $g\geq d/2$ we obtain $h^0(X,L)\leq d/4+1/2$, so we
are done.
If $d_1\geq 0$, by Lemma~\ref{e} we have 
$$
h^0(X,L)\leq  d+1-\min \{d_2, g\}.
$$
If $d_2<g$, we obtain 
$$
h^0(X,L)\leq d+1- d_2=d_1+1\leq d/2+1
$$
  (as  $d_1\leq d/2$); so
we are done. If $d_2>g$, then 
$$h^0(X,L)\leq d+1- g\leq d/2+1
$$ (as $g\geq d/2$).
If $d=0$ and 
 $h^0(L)=1$, by Proposition~\ref{empty} we need to have $\mdeg L\geq 0$.
By Corollary 2.2.5 in \cite{Ctheta} we get  $L=\O_X$.
Finally, suppose $d=2g-2$ and 
let $L$ be balanced, such that $h^0(L)=g$; by 
Serre duality $h^0(\omega_X\otimes L^{-1})=1$. By the previous case and Remark~\ref{balr} (\ref{balrd}).
$\omega_X\otimes L^{-1}=\O_X$, so the proof of Proposition~\ref{cls}  is done.

For  part~\ref{clss} 
let $\hL\in \Pic \hXS$ with $\#S=e\geq 1$
(notation in Subsection~\ref{notBN}).
 We have
$ 
h^0(\hXS,\hL)=h^0(Y_S,M) 
$ (by (\ref{key})),
where $M=\hL_{|Y_S}$ has degree $d-e<d$. 
If $e\leq g-1$ then $Y_S$ has genus at least $1$ so the result follows from (Proposition~\ref{cls}) applied to $Y_S$, which we
can do because $Y_S$ is a binary curve and $M$ is balanced (cf. Definition~\ref{urb}). Otherwise
$Y_S$ has genus $0$ in case $e=g$, or $-1$ if $e=g+1$. In both cases we get
$ 
h^0(Y_S,M)\leq d-g+1\leq d/2+1
$.\end{proof}

Let now $0<d<2g-2$, recall that for a smooth curve $C$, there exists $L\in \Pic^dC$
with $h^0(L)=d/2+1$ if and only if $C$ is hyperelliptic and $L$ is a multiple of the hyperelliptic
class.
The analogous fact holds for binary curves, as we shall see in Theorem~\ref{bcl}.
First we need to define and study hyperelliptic binary curves.

Let $X$ be a binary curve of genus $g\geq 2$. 
$X$ (like all stable curves, cf. \cite{HM})  is called  {\it  hyperelliptic},
if
$X$ lies in the closure, $\Hgb\subset\Mgb$, of the locus, $H_g$, of smooth hyperelliptic curves.
We say that  $X$ is {\it weakly hyperelliptic}
if $W^1_{\md}(X)\neq \emptyset$ for some balanced $\md$ with $|\md|=2$.
If $g\leq 1$ we  say that every binary curve is hyperelliptic (and weakly hyperelliptic), for simplicity.

\begin{remark}\label{whr}
By Proposition~\ref{empty}, $X$ is weakly hyperelliptic if and only if $W^1_{(1,1)}(X)\neq \emptyset$.
\end{remark}

\begin{lemma}
\label{whyp} Let $X$ be a binary curve of genus $g\geq 2$.

\begin{enumerate}[(i)]
\item
\label{whyp1}
$X$ is weakly hyperelliptic if and only if it is hyperelliptic.
\item
\label{whypdef}
If $X$ is hyperelliptic, then   $W^1_{(1,1)}(X)=\{H_X\}$;    $H_X$ will be called the {\emph {hyperelliptic class}}
of $X$.
\item
\label{whyp3}
If $X$ is hyperelliptic, every normalization of $X$ is hyperelliptic.
If $g\geq 4$ and $X$ is not hyperelliptic, there exists a node $n\in \sing$ such that the normalization of $X$ at $n$
is not hyperelliptic.
\end{enumerate}
\end{lemma}
\begin{proof}
Suppose $X$ hyperelliptic, then $\Wudb\neq \emptyset$, by Proposition~\ref{Wp} (\ref{Wp3}).
To show that $X$ is weakly hyperelliptic,
we need to prove  
$W^1_{\md}(X)\neq \emptyset$, for some $\md\in \BXd$.
Pick $[M,S]\in \PXdb$
with $S\neq \emptyset$;
it suffices to show that  $h^0(Y_S,M)\leq 1$.

As $\#S=e\geq 1$ we get $\deg M=2-e\leq 1$. We also know that $\mdeg M$ is balanced, by Definition~\ref{urb}.
By  Proposition~\ref{cliff}, we have 
$$
h^0(Y_S,M)\leq \deg M/2+1\leq 3/2
$$ hence $h^0(Y_S,M)\leq 1$.

Conversely, let $X$ be weakly hyperelliptic. By Remark~\ref{whr} this is equivalent to saying that 
$W^1_{(1,1)}(X)\neq \emptyset$, so $X\in B^1_{g,2}$
(notation in (\ref{Bgrd})). On the other hand, every $X'\in  B^1_{g,2}$  has 
$W^1_{(1,1)}(X')\neq \emptyset$.
Therefore  $B^1_{g,2}=B^1_{g,(1,1)}$; now $B^1_{g,(1,1)}$ is easily seen to be irreducible of dimension $g-2$.

Consider $\Hgb\subset \Mgb$, the locus of   hyperelliptic stable curves.
By the previous part $\Hgb\cap B_g \subset B^1_{g,2}$, hence
\begin{equation}
\label{dimBH}
\dim \Hgb\cap B_g\leq \dim B^1_{g,2}=g-2.
\end{equation}
On the other hand, as $B_g$ is irreducible of codimension $g+1$ in $\Mgb$ 
(cf. Lemma~\ref{Bg}) we have
$$
\dim \Hgb \cap B_g\geq \dim \Hgb - (g+1)=g-2.
$$
Combining this with (\ref{dimBH}) we obtain $\dim \Hgb\cap B_g=g-2=\dim  B^1_{g,2}$.
Since $ B^1_{g,2}$ is irreducible and contains $\Hgb \cap B_g$, we conclude $\Hgb \cap B_g= B^1_{g,2}$, proving
(\ref{whyp1}).

Now, suppose $X$ hyperelliptic, so that $W^1_{(1,1)}(X)\neq \emptyset$.
To prove (\ref{whypdef})
we use induction on $g$:
if $g=2$,  by Proposition~\ref{cliff}, $W^1_{(1,1)}(X)$ contains a unique element: $\omega_X=H_X$.
Now, let $g\geq 3$ and $Y\to X$   the normalization of one node of $X$,
so that $g_Y=g-1\geq 2$.
As $W^1_{(1,1)}(X)\neq \emptyset$ the pull-back map
$$
\rho:W^1_{(1,1)}(X) \la W^1_{(1,1)}(Y);\  \  \  L\mapsto \nu^*L
$$
shows that $Y$ is also weakly hyperelliptic, hence hyperelliptic. 

By induction $W^1_{(1,1)}(Y)=\{H_Y\}$;
by Proposition~\ref{cliff}, $h^0(Y,H_Y)=2$, and $H_Y$ has no base points. Therefore,
by Lemma~\ref{dsl},   $\rho^{-1}(H_Y)$   is a point,
so we are done.

For the final part,
it remains to show that if $X$ is non-hyperelliptic and $g\geq 4$, there exists  $n\in \sing$
such that the normalization   at $n$ is not hyperelliptic.
By contradiction, suppose  this is not the case. Let $Z\to X$ be the normalization of $X$ at two nodes,
$n_1, n_2$,
and call $Y_i$ the normalization of $X$ at $n_i$;
so $Y_i$ is hyperelliptic, for $i=1,2$.
 Therefore $Z$ is hyperelliptic (by the previous part) and has genus at least $2$. Hence $W^1_{(1,1)}(Z)=\{H_Z\}$
and, as $Y_i$ is hyperelliptic,
$$
p_i\sim _{H_Z}q_i,\  \  \  i=1,2,
$$
where $p_i,q_i\in Z$ are the branches over $n_i$.
But then, by Lemma~\ref{dsl}, there exists $L\in \picX{(1,1)}$ which pulls back to $H_Z$ and such that 
$h^0(X,L)=2$. Hence $X$ is weakly hyperelliptic, and hence hyperelliptic
(by the previous part), a contradiction.
\end{proof}

\subsection{Clifford index}
Recall that the {\it Clifford index} of a line bundle $L$ on a   curve $X$ is   
$ 
\Cl  L: =\deg L - 2h^0(L)+2.
$  
Let us define the Clifford index of    $X$:
\begin{equation}
\label{ClX}
\Cl  X :=\min\{\Cl L | L\in \Pic X,\  \mdeg L\in \BXd, \  h^0(L)\geq 2,\  h^1(L)\geq 2 \}.
\end{equation}

For a smooth curve $C$, $\Cl C\geq 0$, and $\Cl C= 0$ if and only if $C$ is hyperelliptic
(Clifford's theorem).
If $C$ is non-hyperelliptic, then 
 $\Cl C= 1$ if and only if $C$ is trigonal or bielliptic  or  a plane quintic
(Mumford's theorem, see   \cite{ACGH} IV (5.2)).

\begin{thm}
\label{bcl}
Let $X$ be a binary curve.
\begin{enumerate}[(I)]
\item
\label{bcl>}
$\Cl X\geq 0$.
\item
\label{bcl0}
 $\Cl X=0$ if and only if $X$ is  hyperelliptic (i.e. weakly hyperelliptic). 
\item
\label{bcl1}
Assume $\Cl X\neq 0$. Then 
 $\Cl X=1$ if and only if $W^1_{\md}(X)\neq \emptyset$ for some balanced $\md$ with  $|\md|=3$.
\end{enumerate}
\end{thm}
Part (\ref{bcl>}) is Proposition~\ref{cliff}.
To prove the rest  we need some auxiliary results.

\begin{lemma}
\label{Wbal} Let $X$ be a binary curve of genus $g\geq 1$;
let $\md =(d_1,d_2) \in \BXd$, with  $0\leq d\leq 2g-2$. 
 Assume $d_1\leq d/2-1$.
Then
$ 
W^{[\frac{d}{2}]}_{\md}(X)=\emptyset.
$ 
\end{lemma}
\begin{proof}
Let $L\in \picX{\md}$ and $l=h^0(X,L)$; it suffices to prove that $l\leq d/2$.

If $d_1<0$ we have
$$ 
l=h^0(C_2, L_2(-C_1\cap C_2))=\max \{0, d_2-g\}.
$$
As $\md$ is balanced, 
by (\ref{BIb}) we have 
$$
d_2-g\leq (d-g+1)/2\leq d/4
$$  
(as $g\geq d/2+1$).
So we are
done.

Let $d_1\geq 0$;   Lemma~\ref{e} yields
$ 
l\leq d+1 - \min \{d_2, g\}.
$ 
If $d_2\geq g$ we get 
$$
l\leq d+1-g\leq d+1-d/2-1=d/2
$$
(again, as  $g\geq d/2+1$). So we are done.
Finally, if $d_2<g$,
$$
l\leq d+1-d_2=d_1+1.
$$
By hypothesis, if $d$ is even,   $d_1\leq d/2-1$, hence   $l\leq d/2$ and we are done.
If  $d$ is odd,   $d_1\leq (d-3)/2$, so that $l\leq (d-1)/2$, so we are done.
\end{proof}
\begin{cor}
\label{bcor}
Let $X$ be a binary curve of genus $g$.
$\Cl X =0$ if and only if there exists
  an integer $h$, \  $1\leq h\leq g-2$, such that $W^h_{(h,h)}(X)\neq \emptyset$.

Assume $\Cl X >0$; then
$\Cl X =1$ if and only if there exists
  an integer $h$, \  $1\leq h\leq g-2$, such that $W^h_{(h,h+1)}(X)\neq \emptyset$.
\end{cor}

\begin{prop}
\label{bva} Let $X$ be a binary curve; its dualizing sheaf,    $\omega_X$, is very ample
if and only if $W^1_{(1,1)}(X)=\emptyset$ (if and only if $X$ is not hyperelliptic).
\end{prop}
\begin{proof}
The part in parentheses follows from Remark~\ref{whr} and Lemma~\ref{whyp}.
Assume $W^1_{(1,1)}(X)=\emptyset$.
We denote   $\dfib:= X\smallsetminus \sing$ the smooth locus of $\dfib$.
For every (not necessarily distinct) $p,q\in X$
we have
\begin{equation}
\label{pq}
h^0(\omega_X(-p-q))=g-3+h^0(X,p+q)=g-2
\end{equation}
($h^0(X, p+q)=1$ by hypothesis and by Lemma~\ref{Wbal}).

Now, for every node $n\in \sing$, denote $\nu:Y\to X$ the normalization at $n$, and $\nu^{-1}(n)=\{r,s\}$;
note that $\omega_Y=\nu^*\omega_X (-r-s) $.
Calling ${\mathcal I}_n$ the ideal sheaf of $n$ in $X$, we have
\begin{equation}
\label{nY}
h^0(X, \omega_X \otimes {\mathcal I}_n)=h^0(Y, \nu^*\omega_X (-r-s))=h^0(Y, \omega_Y)=g-1.
\end{equation}
Formulas (\ref{pq}) and (\ref{nY}) yield that $\omega_X$ is globally generated  and induces a morphism
$\phi:X\to \pr{g-1}$ whose restriction to $\dfib$ is an immersion.
It remains to prove that $\phi$ is injective, and an immersion locally at the singular points of $X$.
Notice that for every nonsingular point $y\in Y$ we have
\begin{equation}
\label{yY}
h^0(Y, \omega_Y(-y))=g-2
\end{equation}
(as $h^0(Y,y)=1$).
Now, for every $p\in \dfib$ and $n\in \sing$ we have, with the same    notation as above
(calling again $p\in Y$  the point  over $p\in X$)
$$
h^0(X, \omega_X (-p)\otimes {\mathcal I}_n)=h^0(Y, \nu^*\omega_X (-p-r-s))=h^0(Y, \omega_Y(-p))=g-2
$$
by (\ref{yY}). Hence $\phi(n)\neq \phi(p)$.
Now let $n_1, n_2\in \sing$, denote $\nu':Y'\to X$ the normalization at $n_1$ and $n_2$, and
$(\nu')^{-1}(n_i)=\{r_i,s_i\}$. We have
$$
h^0(X, \omega_X \otimes {\mathcal I}_{n_1}\otimes {\mathcal I}_{n_2})=
h^0(Y, \nu^*\omega_X (-r_1-s_1-r_2-s_2))=h^0(Y',
\omega_{Y'})=g-2.
$$
Therefore $\phi$ is injective. To show that $\phi$ is
an immersion  at every $n\in \sing$
it suffices to show 
that 
$ 
 H^0(Y, \nu^*\omega_X (-2r-2s))\neq  H^0(Y, \nu^*\omega_X (-2r-s))
$  and that
$ 
 H^0(Y, \nu^*\omega_X (-3r-s))\neq  H^0(Y, \nu^*\omega_X (-2r-s))
$
(notation as above). By (\ref{yY}), 
$$h^0(Y, \nu^*\omega_X (-2r-s))=g-2.$$
On the other hand 
$$
h^0(Y, \nu^*\omega_X (-2r-2s))=h^0(Y,  \omega_Y (-r-s))=g-4+h^0(Y,r+s)=g-3,
$$ 
indeed, if we had $h^0(Y,r+s)=2$ then, by Lemma~\ref{dsl}, $W_{(1,1)}^1(X)$ would be non empty, which is impossible. 
Similarly, $h^0(Y, \nu^*\omega_X (-3r-s))=g-4+h^0(Y,2r)=g-3$ by Proposition~\ref{empty}.
This finishes the first half of the proof.

The opposite implication is easy; let $\omega_X$ be very ample.   By contradiction, let $L\in W^1_{(1,1)}(X)$.
For any $p\in \dfib$ we have $h^0(L(-p))=1$.
So, $L=\O_X(p+q)$ for some $p,q\in \dfib$. Hence $h^0(\omega_X(-p-q))=g-1$,
contradicting the very ampleness of $\omega_X$.
\end{proof}

\begin{lemma}
\label{Wsp} Let $X$ be a binary curve of genus $g\geq 3$ with  $\omega_X$   very ample.
Then
\begin{enumerate}[(i)]
\item
\label{Wsp1}
$W^h_{(h,h)}(X)=\emptyset$ for every $2\leq h\leq g-2$.
\item
\label{Wsp2}
$W^h_{(h,h+1)}(X)=\emptyset$ for every $2\leq h\leq g-4$.
\end{enumerate}
\end{lemma}
\begin{proof}As  $\omega_X$ is very ample, we  identify $X$ with its canonical model
in $\pr{g-1}$, which is a union of
two rational normal curves, $C_1$ and $C_2$   meeting transversally at $g+1$ points.
By contradiction, let $L\in W^r_{\md}(X)$, with $(r,d)$   as in the statement.

We claim that there exists $D\in \Div X$, $D\geq 0$, 
$D$ supported on the smooth locus of $X$, such that $L=\O_X(D)$.
By contradiction, assume   there is a node $n\in \sing$
such that $s(n)=0$ for every $s\in H^0(X,L)$. Denote $\nu:Y\to X$ the normalization of $X$ at $n$,
so that $Y$ is a binary curve of genus $g-1\geq 2$. Set $\nu^{-1}(n)=\{p,q\}$, and $M=\nu ^*L$.
By assumption,  $h^0(Y,M(-p-q))\geq h^0(X,L)$,
therefore  $h^0(Y,M(-p-q))\geq h+1$.
On the other hand, $\mdeg  M(-p-q)=(h-1,h-1)$ is obviously balanced; furthemore
$\mdeg M\geq 0$, hence  Clifford's
inequality  yields
$h^0(Y,M(-p-q))\leq h$, a contradiction. The claim is proved.

Fix  such a $D$, and denote by $\Lambda  \subset \pr{g-1}$ the linear subspace spanned
by
$D$ (if $D$ is  reduced  $\Lambda$ is the ordinary linear span of the points of $D$, 
otherwise $\Lambda$ is the linear span of the appropriate osculating spaces of $X$  
at the points of $\supp D$).
The geometric version of the Riemann-Roch Theorem
 (\cite{ACGH} p. 12) yields
\begin{equation}
\label{GRR}
h^0(X,L)=\deg L - \dim \Lambda .
\end{equation}

In case (\ref{Wsp1}), since $h^0(X,L)\geq h+1$ and  $\deg L=2h$,  we get
\begin{equation}
\label{GRR1}\dim \Lambda\leq h-1.
\end{equation}
If 
$h=1$, then   $\dim \Lambda =0$, which is impossible, as $\Lambda$ is spanned by two distinct points
(as $\mdeg D=(1,1)$). 
So we can assume  $h\geq 2$. We denote
$D=\sum_{i=1}^h(r_i+s_i)$ with $r_i\in C_1$ and $ s_i\in C_2$.
We have  $h^0(X,D -r_1)\geq h+1-1\geq 2$,  hence there exists an effective divisor $D'\neq D$,
with $D\sim D'$, $\supp D'\subset \dfib$, and such that $r_1 $ is in the support of $D'$.
Let $\Lambda ' $ the linear subspace spanned by $D'$ and 
 $\Gamma=<\Lambda ,\Lambda '>$. We have $\dim \Lambda '\leq h-1$ 
and 
$$
\dim \Gamma \leq 2h-1-c
$$
where $c$ is the degree of the greatest common (effective) divisor of $D$ and $D'$; thus $c\geq 1$, by construction.
Now  
 we have, as $r_1\not\in C_2$,
$$
\deg \Gamma\cdot C_2\geq h+h-c+1=2h-c+1,
$$
and this is impossible: $C_2$ is a rational normal curve, so $\Gamma$ cuts on it a divisor of degree at most 
$\dim \Gamma +1=2h-c$.

For part~(\ref{Wsp2}) the method is essentially the same.  
By (\ref{GRR}) 
we have  
$\dim \Lambda\leq h 
$
 and $\Lambda$ is an $(h,h+1)$-secant space of $X$.
Set
$D=\sum_{i=1}^h(r_i+s_i)+s_{h+1}$ with $r_i\in C_1$ and $ s_i\in C_2$.
We have  $h^0(X,D -r_1)\geq 2$, hence there is an effective   $D'\neq D$,
 $D\sim D'$, with $D'-r_1\geq 0$.
With the same notation as above,  $\dim \Lambda '\leq h$ 
 and 
$ 
\dim \Gamma \leq 2h+1-c,
$ 
where $c\geq 1$ was defined above.

Now, 
$ \deg \Gamma\cdot C_2\geq 2h+2-c+1=2h-c+3,
$ 
a contradiction.
\end{proof}

{\it   End of the proof of Theorem~\ref{bcl}.}
Part (\ref{bcl0}). 
By Lemma~\ref{whyp}, $X$ is hyperelliptic if and only if it is weakly hyperelliptic.
If $X$ is weakly hyperelliptic, then $\Cl X=0$.
We prove the converse by showing that if $X$ is not weakly hyperelliptic, then $\Cl X>0$.
By Corollary~\ref{bcor}, it is enough to prove that 
$W^h_{(h,h)}(X)= \emptyset$
for every $h$ with $1\leq h\leq g-2$.

To say that $X$ is not weakly hyperelliptic is to say that $W^1_{(1,1)}(X)=\emptyset$.
By Proposition~\ref{bva}, this implies that $\omega_X$ is very ample.
Lemma~\ref{Wsp} yields $W^h_{(h,h)}(X)= \emptyset$, as wanted.
The proof of part~(\ref{bcl0}) is  complete.

For part (\ref{bcl1}), one direction is obvious.
For the converse, suppose $W^1_{\md}(X)= \emptyset$  for every $\md \in \BXt$ and let us prove that $\Cl X> 1$.
As we are also assuming $\Cl X\neq 0$ we have     $W^1_{(1,1)}(X)=\emptyset$,
hence $\omega_X$ is very ample. Lemma~\ref{Wsp} (\ref{Wsp2}) yields  
$W^h_{(h,h+1)}(X)= \emptyset$ for every $2\leq h\leq g-4$. By Lemma~\ref{Wbal}
it remains to show that $W^1_{(1,2)}(X)$ and $W^{g-3}_{(g-3,g-2)}(X)$ are empty.
The former is empty by assumption; the latter is empty because the former is 
(by Serre duality).   Theorem~\ref{bcl} is proved.
$\blacksquare$
\subsection{Extension of Martens theorem}
\begin{lemma}
\label{HX} Let $X$ be a hyperelliptic binary curve of genus $g\geq 2$,
and $L\in \picX{d}$ be balanced, with $0\leq d\leq 2g-2$.  Then
$\Cl L=0$
if and only if $L=H_X^{\otimes\frac{d}{2}}$
($H_X$ as in Lemma~\ref{whyp}).
\end{lemma}
\begin{proof} 
By the base-point-free-pencil trick
we have
 $h^0(X,H_X^{\otimes\frac{d}{2}})=d/2+1,
 $
so that $\Cl H_X^{\otimes\frac{d}{2}}=0$.
If $g=2$ the statement was proved in Proposition~\ref{cliff}. We continue by induction on $g$.
If $d=2g-2$ then $L=\omega_X$, hence $\omega_X=H_X^{g-1}$.
So we can further assume $d\leq 2g-4$. Let $\md=\mdeg L$, so that $\md\in \BXd$.
By Proposition~\ref{empty} we must have $\md=(d/2,d/2)$; set $r=d/2$.
Let $\nu:Y\to X$ be the normalization of $X$ at  one node,
then
$\nu^*L\in W_{(r,r)}^r(Y)$. Obviously $(r,r)$ is balanced on $Y$.
By induction $W_{(r,r)}^r(Y)=\{H_Y^r\}$, and $h^0(Y,H_Y^r)=r+1$
by Clifford. By Lemma~\ref{dsl},
$ W_{(r,r)}^r(X)$ contains at most one element, hence $L=H_X^r$.
\end{proof}
Martens  Theorem holds for  binary curves, by the following Proposition.
\begin{prop}
\label{martens}
Let $X$ be a binary curve of genus $g\geq 3$. Fix $d,r$ such that $2\leq d\leq g-1$ and $0<2r\leq d$.
Let $\md=(d_1,d_2)\in \BXd$ and assume $r\leq d_i$ for $i=1,2$ (otherwise 
$\dim\Wmdr=\emptyset$, by Prop.~\ref{empty}).

If $X$ is not hyperelliptic,
then $\dim\Wmdr\leq d-2r-1$. 

If $X$ is  hyperelliptic,
then $\dim\Wmdr= d-2r$.
\end{prop}
\begin{proof}
Recall that   if $X$ is hyperelliptic then $W^1_{(1,1)}(X)=\{H_X\}$;
if $X$ is not hyperelliptic, then  $W^1_{(1,1)}(X)$ is empty.
We use induction on $g$.

If $g=3$ then $d=2$ and $r=1$, so the only case to consider is $\md =(1,1)$.
If $X$ is hyperelliptic,   $W^1_{(1,1)}(X)=\{H_X\}$   so it is irreducible of dimension $0$, as
claimed. If $X$ is not hyperelliptic, then $W^1_{(1,1)}(X)=\emptyset$, so we are done.

Let   $g\geq 4$.
If $X$ is not hyperelliptic, by Lemma~\ref{whyp} there exists a node $n\in \sing$ such 
that the normalization $\nu:Y\to X$ of $X$  at $n$ is not hyperelliptic.
Suppose $\Wmdr\neq \emptyset$;
consider the pull-back map
$$
\rho:\Wmdr  \la W^r_{\md}(Y);\  \    \  L\mapsto \nu^*L.
$$
Notice that $\md\in B_d(Y)$; indeed if (say) 
$$
d_1<m(d,g-1)=\frac{d-(g-1)-1}{2}=\frac{d-g}{2}\leq\frac{g-1-g}{2}
$$
(as $d\leq g-1$). So $d_1<0$, hence $\Wmdr=\emptyset$. A contradiction.

If $d\leq g-2=g_Y-1$ we  use induction to get
$\dim  W^r_{\md}(Y)\leq d-2r-1$ and
$\dim  W^{r+1}_{\md}(Y)\leq d-2r-3$. 
Now, suppose  $W^r_{\md}(Y)$
does not have the two points $\nu^{-1}(n)$ as  fixed base points. Then
the fibers of $\rho$ over $W^r_{\md}(Y)\smallsetminus W^{r+1}_{\md}(Y)$
have dimension $0$ (by Lemma~\ref{dsl}), and over $W^{r+1}_{\md}(Y)$ have dimension at most $1$.
Therefore $\dim \Wmdr\leq d-2r-1$. 

If instead $\nu^{-1}(n)$ are  base points of every element of $W^r_{\md}(Y)$,
then, by induction, $\dim W^r_{\md}(Y) \leq (d-2)-2r-1=d-2r-3$ and hence
$\dim \Wmdr\leq d-2r-2$. The case $d\leq g-2$ is settled.

Now let $d=g-1$; then, by Serre duality, $W^r_{\md}(Y)\cong W^{r-1}_{\me}(Y)$ where 
$$\me=(g_Y-1,g_Y-1)-\md\in B_{g_Y-2}(Y),$$
by Remark~\ref{balr} (\ref{balrd}).
Therefore, by induction,
$$
\dim W^r_{\md}(Y)=\dim W^{r-1}_{\me}(Y)\leq g_Y-2-2r+2-1=(g-1)-2r-1
$$
and $\dim  W^{r+1}_{\md}(Y)=  \dim W^{r}_{\me}(Y)\leq (g-1)-2r-3$.
Arguing as before we are done.

Let now $X$ be hyperelliptic.
Then $W_{(r,r)}^r(X)=\{H_X^r\}$ by Lemma~\ref{HX}. Therefore the statement holds if $d=2r$, and
we can assume  $d>2r$.
An induction argument, analogous to the previous one, shows that
$\dim \Wmdr\leq d-2r$
(now $Y$ is  hyperelliptic). To prove that equality holds, pick
 $x_1,\ldots ,x_{d-2r}\in \dfib$
such that $\mdeg H_X^r(\sum x_i)=\md$. It is clear that
$H_X^r(\sum_1^{d-2r} x_i)\in \Wmdr$. Moreover, by Lemma~\ref{abel},
$$
H_X^r(\sum_{i=1}^{d-2r} x_i)\not\cong H_X^r(\sum_{i=1}^{d-2r} x_i'),
$$
 for $x_i$ and $x'_i$ generic.
This shows that $\dim\Wmdr\geq d-2r$, finishing the proof.
\end{proof}
Suppose $d=g-1$, then 
$$\ov{W^0_{g-1,X}}=\Theta(X)
$$ 
where $\Theta(X)$ is the Theta divisor, known to be Cartier and ample (\cite{alex}).
It is thus worth pointing out the following special case of Proposition~\ref{martens}. 
\begin{remark}
\label{theta}
Let $X$ be a binary curve of genus $g\geq 3$. 
For every multidegree $\underline{g-1}\in B_{g-1}(X)$ with $\underline{g-1} >0$ we have
\begin{displaymath}
\dim W_{\underline{g-1}}^1(X)=\left\{ \begin{array}{l}
g-3 \   \text{ if }\  X  \text{ is hyperelliptic } \\
g-4 \   \text{ otherwise.} \\
\end{array}\right.
\end{displaymath}
If $X$ is  an irreducible curve
the same holds (\cite{Ctheta} Thm. 5.2.4).
\end{remark}

\section{Dimension of Brill-Noether varieties.}
\label{BNsec}

The Brill-Noether number  $\BN$ is defined as follows 
\begin{equation}
\label{rho}
\BN=g-(r+1)(g-d+r)=(r+1)d-rg-(r+1)r.
\end{equation}
By the famous Brill-Noether  theorem,  
$\dim W^r_d(C)=\BN$
for a general smooth curve $C$.
The proof of this theorem has an   interesting history, as many mathematicians have contributed to it: 
Arbarello, Cornalba, Eisenbud, Gieseker, Griffiths, Harris, Kempf, Kleiman, Laksov, Lazarsfeld, Martens,
among others.
We refer to Chapter 5 of \cite{ACGH} for details and references.
 
The  goal of this section is to prove it for   binary
curves, assuming $r\leq 2$. More precisely, we shall prove that for a general binary curve $X$ of genus $g$ and every
balanced multidegree $\md\in \Bgd$ we have  $\dim \Wmdr\leq \BN$, with equality holding for certain $\md$.
As a by-product we have a new proof for smooth curves.
More generally,  
 Theorem~\ref{BNt}  implies that
 the Brill-Noether  theorem holds
 on every stratum of $\Mgb$  containing $B_g$ in its closure.

\begin{thm}
\label{BNt}
Let $X$ be a general binary curve of genus $g$; fix $r\leq 2$ and $d\in \Z$; let $\md \in \BXd$.
Then 
\begin{enumerate}[(i)]
\item\label{BNi}
$\dim \Wmdr \leq  \BN$  and equality holds for some $\md$.
\item
\label{BNii}
$\dim \Wb= \BN$.
\end{enumerate} 
\end{thm}
\begin{proof}
We have $\dim \Wb\geq \BN$
(by Theorem V (1.1) in \cite{ACGH}, which is independently due to 
 \cite{kempf} or \cite{KL1});
also,    $\rho^r_{d-1}(g-1)=\BN-1$.
Therefore, by (\ref{Wstr}), part (\ref{BNii}) follows from part (\ref{BNi}).
So  it   suffices to prove $\dim \Wmdr\leq \BN$.

If $d\geq r+g$ then $\BN\geq g$, so the statement is trivial. We shall thus assume
$d\leq r+g-1$. 
If $d\leq 0$, then $\Wmdr=\emptyset$ (by Proposition~\ref{empty}), unless $\md=(0,0)$, in which case
$\Wmdr=\{\O_X\}=W^0_{\md}(X)$
(by Corollary 2.2.5 in \cite{Ctheta}), so the theorem holds.

We can thus use induction on $d$.
We set $d_1\leq d_2$.
By Lemma~\ref{empty}, $\Wmdr=\emptyset$ if $d_1\leq r-1$; therefore we can assume $d_i\geq r$.

We begin with  $r=0$, in which case   a more precise result holds.
Recall   that we called  $\Amd\subset  W_{\md}(X)$
the closure of the image of the $\md$-th Abel map; see Lemma~\ref{abel}.
\begin{prop}
\label{BN0} Let $X$ be any binary curve of genus $g$,  and $d\leq g-1$.
For every $\md \in \BXd$ with $\md \geq 0$,  $\dim  W_{\md}(X) = d$.
Moreover $\Wmd$ has a unique irreducible component of dimension $d$, namely $\Amd$.
\end{prop}
\begin{proof}
Suppose $d=g-1$. By Theorem~3.1.2 of \cite{Ctheta},  
if $\md$ is strictly balanced
the proposition holds.
For a binary curve, the only  balanced,   non strictly balanced,  multidegree
is $(-1,g)$, which is
 ruled out by hypothesis.
The case $d=g-1$ is done.

We continue by induction on $g-d$.
Let $d\leq g-2$ and consider the normalization $\nu:Y\to X$ of $X$ at one node, $n$; set
$\nu^{-1}(n)=\{p,q\}$.
Thus $Y$ is a binary curve of genus $g-1$.
Consider  
$$
\rho:\Wmd \la W_{\md}(Y);\  \  \   L\mapsto \nu^*L.
$$
If $\md$ is not  balanced for $Y$,  we may assume (up to switching $C_1$ and $C_2$)
$$
d_1< \frac{d-(g-1)-1}{2}\leq \frac{g-2-g}{2}=-1;
$$ impossible. So   $\md$  is balanced for $Y$.
Thus, by induction,   $W_{\md}(Y)$ has a unique component of dimension $d$,
namely $A_{\md}(Y)$. 

Call $B\subset A_{\md}(Y)$ the locus of $M$ such that
$h^0(M)=h^0(M(-p))=h^0(M(-q))=1$. 
Then one easily checks that  $\dim B\leq d-2$, hence $\dim \rho^{-1}(B)\leq d-1$
(since the fibers of $\rho$ have dimension at most $1$, of course).

By Lemma~\ref{abel}, there exists a dense open subset
$U\subset  A_{\md}(Y)\smallsetminus B$ such that
for every $M\in U$ we have $h^0(Y,M)=1$.
By Lemma~\ref{dsl} the fibers of $\rho$ over such and $M$ is a unique point.
Therefore $\rho^{-1}(U)$ is irreducible of dimension $d$

Now, by Proposition~\ref{martens}, $\dim W^1_{\md}(Y)\leq d-2$, therefore any other component
(if it exists) of $\Wmd$ has dimension at most $d-1$.
This proves that $\Wmd$ has a unique component, $W$, of dimension $d$;
by Lemma~\ref{abel}, $W=\Amd$.
\end{proof}
We point out a simple consequence.
\begin{cor}
\label{cor0} Let $X$ be any binary curve of genus $g$, and $d=g+r-1$.
Then $\dim  W_{\md}(X) = \BN=g-r-1$
for every $\md \in \BXd$ with $d_i \geq r$ for $i=1,2$. 
Moreover $\Wmd$ has a unique irreducible component, $W$, of dimension $\BN$, and 
for the general  $L\in W$ we have
$h^0(X,L)=r+1$.
\end{cor}

\begin{proof}
Set 
$ 
(d_1',d_2')=\md'=\mdeg \omega_X-\md=(g -1-d_1, g-1-d_2).
$ 
We have 
$ 
d_2'= g-1-d+d_1\geq g-1-d+r=0
$; similarly $d_1'\geq 0$, hence $\md '\geq 0$.
Also, $\md'$ is balanced, because $\md$ is (by Remark~\ref{balr} (\ref{balrd}).
By Serre duality,  $\Wmdr \cong W^0_{\md '}(X) $, and
 the corollary follows from Proposition~\ref{BN0} and Lemma~\ref{abel}.
\end{proof}

We now go on with the proof of the theorem.

\begin{claim}
\label{cgg}
Assume $r\geq 1$. Let  $W$ be an irreducible component of $\Wmdr$ having maximal dimension.  Then    the general  $L\in W$ is globally generated. 
\end{claim}
We have $\dim W\geq \BN$.
By contradiction, suppose that every section of $L$ vanishes at $p\in X$.
If $X$ is smooth at $p$, then   $L(-p)$ is a line bundle of multidegree $\md'=(d_1-1,d_2)$
(say). We claim $\md'$ is balanced. Indeed
if $\md'\not\in B_{d-1}(X)$ we must have (since $m(d,g)>m(d-1,g)$)
$$
d_1-1<m(d-1,g)=\frac{d-g-2}{2}\leq \frac{r-3}{2} 
$$
(using $d\leq g+r-1$). Therefore
$ 
d_1< (r-1)/2,
$ 
hence
$ 
d_1\leq r-1,
$ 
 a contradiction.
So, $\md'$ is balanced;  induction  yields
$$
\dim W^r_{\md'}(X)\leq \rho^r_{d-1}(g)=\BN -r-1.
$$
Therefore the set of line bundles in $\Wmdr$ admitting a base point has dimension at most
$\rho^r_{d-1}(g)+1=\BN-r\leq \BN$
(consider the rational  map $\Wmdr\times X \da \picX{\md'}$ mapping $(L,p)$ to $L(-p)$).
So,   $\dim W< \BN$, a contradiction.

Now assume that every section of $L$ vanishes at a node $n$ of $X$. 
Since $X$ has finitely many nodes, we may assume that the node $n$ is the same for the general $L$.
Let $\nu:Y\to X$ be the normalization of $X$ at $n$, so that $Y$ is a binary curve of genus $g-1$.
Denote $\nu^{-1}(n)=\{p,q\}$. Then 
$\nu^*L(-p-q)\in W^r_{d_1-1, d_2-1}(Y)$.
We claim that $(d_1-1, d_2-1)$ is balanced on $Y$. If that were not the case, then (say)
$$
d_1-1<m(d-2,g-1)=\frac{d-2-g+1-1}{2}=\frac{d-2-g}{2}\leq \frac{r-3}{2}.
$$
Therefore (as before)
$ 
d_1< (r-1)/2,
$ 
hence $d_1\leq r-1$; a contradiction.
As  $(d_1-1, d_2-1)$ is balanced  we get (by induction)
$$
\dim W^r_{(d_1-1, d_2-1)}(Y)\leq \rho^r_{d-2}(g-1)=g-1-(r+1)(g+r-d+1)=\BN - r-2.
$$
Now consider the map 
$$\Wmdr \to \Pic^{(d_1-1, d_2-1)}Y; \  \  \  L\mapsto \nu^*L(-p-q).
$$
Its fibers have dimension at most $1$, of course.
The restriction of the above map to $W$ maps  the general element of $W$ in $W^r_{d_1-1, d_2-1}(Y)$, hence
$$
\dim W\leq \dim W^r_{d_1-1, d_2-1}(Y)+1\leq \BN - r-2+1<\BN,
$$
which is impossible.The claim is proved. 

\

\noindent
{\it Proof  of  Theorem~\ref{BNt} for $r=1.$}

 For $i=1,2$
consider the moduli spaces
$M_{0,g+1}(\pr{1},d_i)$; 
 there are natural maps
\begin{equation}
\begin{array}{lccr}
\epsilon_i:&M_{0,g+1}(\pr{1},d_i)&\la &(\pr{1})^{g+1} \  \  \  \  \  \\
&(\phi_i, (p_1,\ldots, p_{g+1})) &\mapsto &(\phi_i(p_1),\ldots, \phi_i (p_{g+1})).
\end{array}\end{equation}
We thus obtain the cartesian diagram:
$$\begin{CD}
V:=M_{0,g+1}(\pr{1},d_1)\times_{(\pr{1})^{g+1}} M_{0,g+1}(\pr{1},d_2) @>\pi_1>> M_{0,g+1}(\pr{1},d_1)\\
@V\pi_2VV @VV\epsilon_1V\\
 M_{0,g+1}(\pr{1},d_2) @>\epsilon_2>> (\pr{1})^{g+1}.\\
\end{CD}$$
As $\epsilon_1$ and $\epsilon_2$ are dominant, we get
\begin{eqnarray}
\nonumber\\
{\dim V =\dim M_{0,g+1}(\pr{1},d_1)+\dim M_{0,g+1}(\pr{1},d_2) -\dim (\pr{1})^{g+1}=}
\nonumber\\
{ 2(d_1+d_2)-4+2(g+1)-g-1=2d+g-3=\rho^1_d(g)+2g-1}
\end{eqnarray}

\

\noindent ($\rho^1_d(g)=2d-g-2$).
Moreover, $V$ is irreducible, as every scheme in the above diagram is so.
Now, 
$V$ has a natural $PGL(2)$-invariant map to $B_g$
$$
\alpha_{\md}: V\stackrel{\psi}{\la} M_{0, g+1}\times M_{0, g+1}\stackrel{\gamma_g}{\la} B_g 
$$
(where $\psi$   forgets the maps to $\pr{1}$).
By Claim~\ref{cgg}, $\alpha_{\md}$ dominates $B^1_{g,\md}$. Furthermore,
for every curve $X\in \alpha_{\md}(V)$ there is a natural, $PGL(2)$-invariant map
\begin{equation}\label{dddm}
 \alpha_{\md}^{-1}(X)\la W^1_{\md}(X);
\end{equation}
this, together with Claim~\ref{cgg}, yields
\begin{equation}\label{ddd}
\dim W^1_{\md}(X)\leq \dim  \alpha_{\md}^{-1}(X)-3.
\end{equation}
Also, $B^1_{g,\md}$ is irreducible
 and 
\begin{equation}\label{B1d}
\dim B^1_{g,\md}\leq  \min \{\dim V-3, \dim B_g\}=\min \{\rho^1_d(g)+2g-4, \dim B_g\}.
\end{equation}
Recall that $\dim B_g=2g-4$.
If $\alpha_{\md}$   dominates $B_g$,  (\ref{ddd}) yields,  for $X$ general,
  $$
  \dim W^1_{\md}(X)\leq \dim V -\dim B_g-3=\rho^1_d(g).
  $$ 
On the other hand if $\alpha_{\md}$ does not dominate $B_g$, $W^1_{\md}(X)$ is  empty.

If $\rho^1_d(g)<0$, by (\ref{B1d}) 
$\alpha_{\md}$ is not dominant, hence
 $W^1_{\md}(X)=\emptyset$ for $X$ general in $B_g$.

If  $\rho^1_d(g)\geq 0$ then, by \cite{KL1} or \cite{kempf}, there exists a $\md$ such that
$\alpha_{\md}$ dominates $B_g$. By what we said above, for every such $\md$, 
  $\dim W^1_{\md}(X) \leq \rho^1_d(g)$ for the general binary curve $X$.
The proof for $r=1$ is complete.

\

\noindent
{\it Proof of  Theorem~\ref{BNt} for $r=2.$}

By Proposition~\ref{cliff} we can  assume $g\geq 3$.
Define
$ 
J\subset M_{0,g+1}(\pr{2},d_1)\times M_{0,g+1}(\pr{2},d_2)
$  as follows
$$
J=\{((\phi_1; p_1,\ldots, p_{g+1}); (\phi_2; q_1,\ldots, q_{g+1}))|  \  \phi_1(p_i)= \phi_2(q_i)\,\  \forall
i=1,\ldots, g+1\}.
$$

Consider the  map
$ 
\Psi:J\stackrel{}{\la}M_{0}(\pr{2},d_1),
$ 
where $\Psi$ is the projection to the first factor
composed with the map forgetting $(p_1,\ldots,
p_{g+1})$.
 
Pick $\phi_1\in M_{0}(\pr{2},d_1)$.  
For every $\phi_2\in M_{0}(\pr{2},d_2)$,  either
$\im\phi_2\cap \im \phi_1$ is a  finite set,
or $\im\phi_1\subseteq\im \phi_2$
(recall that $d_1\leq d_2$); this second case occurs     only if
$d_2=cd_1$  for some $c\geq 1$.
We   partition $J=J_a\cup J_b$, where
  $J_a$ parametrizes points such that $\im\phi_1\not\subset\im \phi_2$, and $J_b=J\smallsetminus J_a$.
So, $J_b=\emptyset$ if and only if $d_1$ does not divide $d_2$, and $J_a= \emptyset$ if and only if $d_1d_2< g+1$.

Assume  $d_1d_2\geq g+1$.  The restriction of $\Psi$ to $J_a$ is dominant and
$\im \phi_1\cap \im \phi_2$  is made of $d_1d_2$ distinct points,
for  $\phi_1$ and $\phi_2$ general. Hence there are finitely many choices
for the $g+1$ marked points, 
$(p_1,\ldots, p_{g+1})$ and $(q_1,\ldots, q_{g+1})$. We conclude
\begin{equation}
\label{Ia}
\dim J_a=\dim M_{0}(\pr{2},d_1)+\dim M_{0}(\pr{2},d_2)=
3d-2 
\end{equation}
(cf. (\ref{maps})).
On the other hand, if $d_2=cd_1$, the fiber of $J_b$ over $\phi_1$
is the set of all $(\phi_2;q_1,\ldots, q_{g+1})$ such that
$\phi_2=\psi\circ \phi_1$ with $\psi\in M_{0}(\pr{1},c).$
Hence
\begin{equation}
\label{Ib}
\dim J_b\leq \dim M_{0}(\pr{2},d_1)+\dim M_{0}(\pr{1},c)+g+1=3d_1 +2c-2+g.
\end{equation}

Now consider $\Bgtmd\subset \Bgtd\subset \Mgb$.
 $J$ has a natural map, $\beta_{\md}$, to
$\Bgtmd$, obtained by restricting to $J$ the composition of the  forgetful map
(disregarding the maps),
with the  map $\gamma_g$
$$
\beta_{\md}: J\la M_{0, g+1}\times M_{0, g+1}\stackrel{\gamma_g}{\la} B_g.
$$

We claim that   $\beta_{\md}(J_b)$ is never dense in $B_g$ (also when $d_1d_2<g+1$).
Notice that
the restriction of $\beta_{\md}$ to $J_b$ forgets both $\phi_1$ and $\psi\in M_{0}(\pr{1},c)$, and it is invariant
with respect to the
$PGL(2)$ diagonal action on $J_b$.
Therefore by (\ref{Ib}), and recalling that $g>2$, we have
$$
\dim\beta_{\md}(J_b)\leq \dim J_b-(3d_1-1)-(2c-2)-3\leq g-2<\dim B_g.
$$

On the other hand   $\beta_{\md}$ restricted to
$J_a$ is  $PGL(3)$-invariant, hence, by
(\ref{Ia}),
$$\dim\beta_{\md}(J_a)\leq \dim J_a-\dim PGL(3)=3d-10=2g-4+\rho^2_d(g)=\dim B_g+\rho^2_d(g)$$ 
(as $\rho^2_d(g)=3d-2g-6$). 

Now we argue  as for $r=1$; note that $\beta_{\md}$ dominates $\Bgtmd$. 
If $\beta_{\md}(J_a)$ is not dense in $B_g$, then, by what we said,  $W^2_{\md}(X)$ is empty for $X\in B_g$ general.
By the above inequality,
 this will  always be the case if $\rho^2_d(g)<0$.  

 As observed at the beginning of the proof, if  $\rho^2_d(g)\geq 0$,  
then $\beta_{\md}$ dominates $B_g$ for some $\md$. For  such $\md$ 
we derive $\dim W^2_{\md}(X) =\rho^2_d(g)$ for the general binary curve $X$.

It remains to handle the case $d_1d_2<g+1$, when $J_a$  is empty.
We proved above that $\dim\Bgtmd\leq \dim \beta(J_b)\leq g-2<2g-4$,
hence $\Wmdr=\emptyset$  for $X$ general in $B_g$.

Theorem~\ref{BNt}  is proved. \end{proof}
\section*{Acknowledgements}
For the writing of this paper I benefitted from several enlightening conversations with Edoardo Sernesi, to whom  I am
grateful. I also   thank Silvia Brannetti for several useful remarks.

\end{document}